\documentclass{sigcomm-alternate}
\usepackage{amsfonts,latexsym,amsmath,amstext,amssymb,verbatim,epsfig,psfrag}
\usepackage{colordvi,pstcol}
\usepackage{graphicx,psboxit}
\usepackage{psfrag}

\begin{document}
\title{Revisiting the D-iteration method: from theoretical to practical computation cost}

\numberofauthors{1}
\author{
   \alignauthor Dohy Hong\\
   \affaddr{Alcatel-Lucent Bell Labs}\\
   \affaddr{Route de Villejust}\\
   \affaddr{91620 Nozay, France}\\
   \email{\normalsize dohy.hong@alcatel-lucent.com}
}

\date{\today}
\maketitle

\begin{abstract}
In this paper, we revisit the D-iteration algorithm in order to better explain  its connection to the Gauss-Seidel method
and different performance results that were observed. In particular, we study here the practical computation cost based on the execution runtime compared to the theoretical number of iterations. We also propose an exact formula of the error for PageRank class of equations.
\end{abstract}
\begin{psfrags}
\section{Introduction}\label{sec:intro}
In this paper, we assume that the readers are already familiar with the idea
of the fluid diffusion associated to the D-iteration \cite{d-algo}
to solve the equation:
$$
X = P.X + B
$$
and its application to PageRank equation \cite{dohy}.

For the general description of alternative or existing iteration methods, one
may refer to \cite{Golub1996, Saad}.

In Section \ref{sec:gs}, we explain the exact connection between the D-iteration
and the Gauss-Seidel iteration. Section \ref{sec:error} presents the formula for
the exact error (distance to the limit with $L_1$ norm) for PageRank type equations.
Section \ref{sec:cs} presents the analysis of the computation cost.

\section{Connection to Gauss-Seidel iteration}\label{sec:gs}

We recall the equation on the $H_n$ history vector associated to
the D-iteration method:
\begin{eqnarray}\label{eq:h}
H_{n} &=& \left(I_d - J_{i_n}(I_d - P)\right)H_{n-1} + J_{i_n} F_0
\end{eqnarray}
where $I_d$ is the identity matrix, $J_k$ a matrix with all entries equal to zero except for
the $k$-th diagonal term: $(J_k)_{kk} = 1$, $F_0$ the initial condition vector (equal to $B$)
and $i_n$ the $n$-th choice of node for the diffusion.
The choice of the sequence $I = \{i_1, i_2, ..., i_n,...\}$ with $i_n \in \{1,..,N\}$
is the main optimization factor of the D-iteration.

The above equation \eqref{eq:h} is in fact equivalent to:
\begin{eqnarray*}
\mbox{If } i\neq i_n &:& (H_n)_i = (H_{n-1})_i\\
\mbox{If } i= i_n &:&(H_{n})_{i_n} = L_{i_n}(P)H_{n-1} + (B)_{i_n}
\end{eqnarray*}
which is exactly the Gauss-Seidel iteration equation if we apply
the diagonal term elimination (division by $1/(1-p_{ii})$).
This means that the history vector $H_n$ we obtain with the D-iteration
is exactly the same results than the Gauss-Seidel's result when applying
the same sequence $I$. In particular, it means that one can apply
the Gauss-Seidel method for any infinite sequence of $I$ and the limit
is not modified.
However, the main difference is that with the D-iteration, we don't use the
equation \eqref{eq:h}. Instead, we use the
column vector $C_i(P)$ to update ($F_n$, $H_n$): the advantage of introducing and working with $F_n$
is that (cf. pseudo-code \cite{d-algo}):
\begin{itemize}
\item we know exactly in advance the amount of fluid on which the diffusion is
applied (so the consequence in advance: how much remains and how much disappears): 
if from $H_n$ (and line vector application), 
we want to optimize the way the sequence $I$ is built, it is not obvious
(and this explains why up to now only the cyclic iteration is done) 
and it would in fact require the information of $F_n$;
\item the computation cost is reduced: this will be illustrated below. The main reason
is that when the diffusion is applied, each computation is {\em useful} (each diffusion
adds fluid effectively to its children nodes), whereas with
the line vector application on $H_n$, we may have a lot of redundant computation.
To understand this last idea, assume that $x$\% of $N$ are constant or almost constant.
With the line vector $L_i(P)$ application, there will be $x$\% of operations that will
be repeated and that could be avoided if the diffusion approach is applied (cf. results
comparison of Section \ref{sec:calif}).
\end{itemize}

\section{Improving the error estimate}\label{sec:error}

It has been shown in \cite{dohy} that $r/(1-d)$ is an exact distance
(for $L_1$ norm) to the limit when $P$ has all columns summing-up tp $d$.
In case of the PageRank equation, we may have zero-column vector in $P$
(if we don't do the $P$ completion operation cf. \cite{deep}). Indeed if
zero-column vector of $P$ (corresponding to dangling nodes) is to be
completed by $d/N$, any iteration scheme would do useless computations.
When, we are working on the $P$ matrix without completion, the limit we
obtain need to be renormalized (by a constant multiplication for diffusion approach
or by constant addition for power iteration).

To take into account this effect precisely, we count the total amount
of fluid that left the system when a diffusion is applied on a dangling node:
we call this quantity $e_n$ (at step $n$ of the D-iteration).
This quantity should have been put in the system by adding $e_n\times d /N$ on
each node, which means that the initial fluid should have been $(1-d-d e_n)/N$
instead of $(1-d)/N$. But then the fluid $d e_n/N$ would have produced
after $n$ steps $(d\times e_n/(1-d))^2/N$ that disappears by dangling nodes, etc.
Applying the argument recursively, the correction that is required on the residual fluid 
$r_n$ (equal to $|F_n|$) is
to replace the initial condition $r_0 = 1-d$ by:
$$
(1-d) + d e_n + d e_n\frac{d e_n}{1-d} + d e_n\left(\frac{d e_n}{1-d}\right)^2 + ...\\
= \frac{(1-d)^2}{1-d - d e_n}.
$$

And $H_n$ need to be renormalized (multiplication) by $(1-d)/(1-d-d e_n)$ so that
the exact $L_1$ distance $|H_{\infty}-H_n|$ is equal to:
$$
|H_{\infty}-(1-d)/(1-d-d e_n) H_n| = r_n/(1-d-d e_n).
$$
Below, in the D-iteration approach, we updated $e_n$ by:
\begin{eqnarray*}
e_n &+=& (F_n)_{i_n}
\end{eqnarray*}
if $i_n$ is a dangling node. 

\section{Analysis of the computation cost}\label{sec:cs}

\subsection{Web graph dataset}
For the evaluation purpose, we used the 
web graph imported from the dataset \verb+uk-2007-05@1000000+
(available on \cite{webgraphit}) which has
41,247,159 links on 1,000,000 nodes.

Below we vary $N$ from $10^3$ to $10^6$ extracting from the dataset the
information on the first $N$ nodes.
Few graph properties are summarized in Table \ref{tab:1}:
\begin{itemize}
\item L: number of non-null entries (links) of $P$;
\item D: number of dangling nodes (0 out-degree nodes);
\item E: number of 0 in-degree nodes: the 0 in-degree nodes are defined recursively:
  a node $i$, having incoming links from nodes that are all 0 in-degree nodes, is
  also a 0 in-degree node; from the diffusion point of view, those nodes are those
  who converged exactly in finite steps;
\item O: number of loop nodes ($p_{ii} \neq 0$);
\item $\max_{in} = \max_i \#in_i$ (maximum in-degree, the in-degree of $i$ is the number of
  non-null entries of the $i$-th line vector of $P$);
\item $\max_{out} = \max_i \#out_i$ (maximum out-degree, the out-degree of $i$ is the number of
  non-null entries of the $i$-th column vector of $P$).
\end{itemize}

\begin{table}
\begin{center}
\begin{tabular}{|l|cccccc|}
\hline
N        & L/N  & D/N   & E/N   & O/N & $\max_{in}$ & $\max_{out}$\\
\hline
$10^3$   & 12.9 & 0.041 & 0.032 & 0.236 & 716   & 130\\
$10^4$   & 12.5 & 0.008 & 0.145 & 0.114 & 7982  & 751\\
$10^5$   & 31.4 & 0.027 & 0.016 & 0.175 & 34764 & 3782\\
$10^6$   & 41.2 & 0.046 & 0     & 0.204 & 403441& 4655\\
\hline
\end{tabular}\caption{Extracted graph: $N=10^3$ to $10^6$.}\label{tab:1}
\end{center}
\end{table}

\begin{figure}[htbp]
\centering
\includegraphics[angle=-90,width=\linewidth]{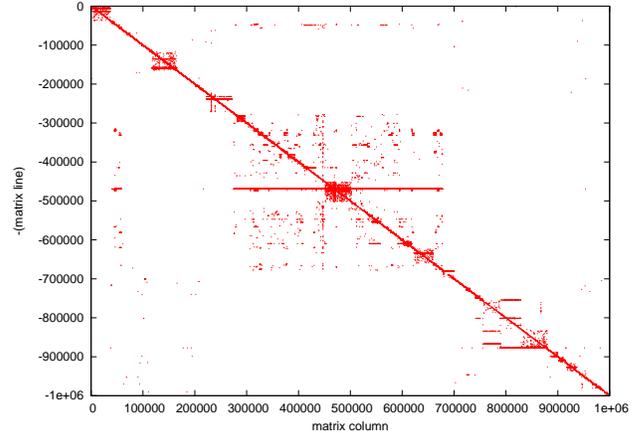}
\caption{P matrix associated to $uk-2007-05@1000000$: random sampling of 100000 links.}
\label{fig:uk}
\end{figure}

\subsection{Programming environment}
For the evaluation of the computation cost, we used $C++$ codes ($g++-3.3$)
on a Linux (Ubuntu) machine:

\verb+Intel(R) Core(TM)2 CPU, U7600, 1.20GHz+, cache size 2048 KB.

The runtime has been measured based on the library {\em time.h} with the function
$clock()$ (with a precision of 10 ms). The runtime below measures the computation time
from the time we start iterations (time to build the iterators are excluded).

\subsection{First comparison tests}
The first algorithms that we evaluated in this section are:
\begin{itemize}
\item PI: Power iteration (equivalent to Jacobi iteration);
\item GS: Gauss-Seidel iteration (cyclic sequence);
\item GS': Gauss-Seidel iteration (cyclic sequence): keeping diagonal terms;
\item DI-CYC: D-iteration with cyclic sequence (a node $i$ is selected, if $(F_n)_i>0$);
\item DI-MAX: D-iteration with $i_n = \arg\max_i (F_{n-1})_i$ by threshold (by threshold means
  that we apply the diffusion to all nodes above the threshold value; when there is no such
  node, we decrease multiplicatively the threshold, by default by 1.2, which we call the
  decrement factor);
\item DI-OP: D-iteration with $i_n = \arg\max_i (F_{n-1})_i/(H_{n-1})_i$ by threshold,
  diffusions are applied first on the 0 in-degree nodes (recursively, such that we need to
  apply the diffusion exactly once to each of those nodes).
\end{itemize}

Because, we are currently limited by the memory size on a single PC, we introduced
an arbitrary function each time a non zero entry of $P$ need to be used by introducing
a finite iteration of 
\begin{verbatim}
void ArbitraryFunction(int m){
  double x = 1.0;
  for (int i=0; i < m; i++){
    x := a * x + b; 
  }
}
\end{verbatim}
($m$ is the number of times we iterates): 
this would correspond exactly to the reality if the operation
$p_{ij}\times x_j$ is to be replaced by an operator $f_{ij}(x_j)$ whose computation cost
is exactly the finite iteration we introduced arbitrarily.

As one can observe in the results of Tables \ref{tab:compa} and \ref{tab:compa2},
the main improvements are brought by:
\begin{itemize}
\item coordinate level update with the iteration GS and DI-CYC (against vector level update of PI);
\item a better choice of the nodes for the diffusion process:
  we see a significant jump with the very basic solution DI-MAX;
\item then with DI-OP (we observed similar improvements with all variants of the idea
  of the argmax of weighted fluid).
\end{itemize}

The impact of the diagonal term elimination and the impact of the redundant computation
of line-vector operations in GS is limited here.

We see also that the computation cost estimate with the number of iterations
is a very good approximation when the matrix product operations (the multiplications)
are the dominant component of the computation run time cost (when $m$ is introduced).

When we set $m=0$ (real runtime), we see that the real speed-up gain may be in fact
much more important than those estimates (for large $N$). Also, we can notice a surprising efficiency of
DI-CYC: if the main speed-up gain is from DI-MAX for the number of iterations,
the maim improvement is brought by DI-CYC for the runtime.

\begin{table}
\begin{center}
\begin{tabular}{|l|ccccc|}
\hline
        & PI  & GS'   & DI-CYC  & DI-MAX  & DI-OP\\
\hline
\hline
$N=10^3$ & & & & &\\
 nb iter & 28  & 22    & 20.8    & 14.4    & 12.8 \\
 speed-up& 1.0 & 1.3   & 1.4     & 1.9     & 2.2 \\
\hline
 $m=10^6$&  & & & & \\
 time (s)& 303 & 238   & 225     & 157     & 138 \\
 speed-up& 1.0 & 1.3   & 1.4     & 1.9     & 2.2 \\
\hline
 $m=0$&  & & & & \\
 time (s)& 0.02 & 0.02   & 0.02     & 0.02     & 0.02 \\
 speed-up& 1.0 & 1.0   & 1.0     & 1.0     & 1.0 \\
\hline
\hline
$N=10^4$ & & & & & \\
 nb iter & 43  & 38    & 34.9    & 17.2    & 13.8\\
 speed-up& 1.0 & 1.1   & 1.2     & 2.5     & 3.1 \\
\hline
 $m=10^5$&  & & & & \\
 time (s)& 453 & 401   & 367     & 181     & 146 \\
 speed-up& 1.0 & 1.1   & 1.2     & 2.5     & 3.1 \\
\hline
 $m=0$&  & & & & \\
 time (s)& 0.29 & 0.28   & 0.13     & 0.23     & 0.11 \\
 speed-up& 1.0 & 1.0   & 2.2     & 1.3     & 2.6 \\
\hline
\hline
$N=10^5$ & & & & & \\
 nb iter & 52  & 43    & 42.9    & 21.1    & 16.2\\
 speed-up& 1.0 & 1.2   & 1.2     & 2.5     & 3.2 \\
\hline
 $m=10^4$&  & & & & \\
 time (s)& 1386& 1146  & 1131    & 560     & 431\\
 speed-up& 1.0 & 1.2   & 1.2     & 2.5     & 3.2 \\
\hline
 $m=0$&  & & & & \\
 time (s)& 9.0 & 8.0  & 2.7     & 3.6     & 2.1\\
 speed-up& 1.0 & 1.1   & 3.3     & 2.5     & 4.3 \\
\hline
\hline
$N=10^6$ & & & & & \\
 nb iter & 66  & 57    & 57      & 18.3    & 15.7 \\
 speed-up& 1.0 & 1.2   & 1.2     & 3.6     & 4.2 \\
\hline
 $m=10^3$&  & & & & \\
 time (s)&2639 & 2275  & 2051    & 693     & 586 \\
 speed-up& 1.0 & 1.2   & 1.3     & 3.8     & 4.5 \\
\hline
 $m=0$&  & & & & \\
 time (s)& 223 & 198   & 45      & 49      & 27\\
 speed-up& 1.0 & 1.1   & 5.0     & 4.6     & 8.3 \\
\hline
\end{tabular}\caption{Comparison of the runtime for a target error of $1/N$. $m=1000000, 100000, 10000, 1000, 0$. GS is computed here without diagonal terms elimination. speed-up: gain factor w.r.t. PI.}\label{tab:compa}
\end{center}
\end{table}

\begin{table}
\begin{center}
\begin{tabular}{|l|ccccc|}
\hline
        & PI  & GS   & DI-CYC & DI-MAX & DI-OP\\
\hline
\hline
$N=10^3$ & & & & & \\
 nb iter & 28  & 18.7 & 17.5   & 13.3   & 11.3\\
 speed-up& 1.0 & 1.5  & 1.6    & 2.1    & 2.5 \\
\hline
 $m=10^6$&  & & & & \\
 time (s)& 303 & 202  & 190    & 144    & 123 \\
 speed-up& 1.0 & 1.5  & 1.6    & 2.1    & 2.5 \\
\hline
 $m=0$&  & & & & \\
 time (s)& 0.02 & 0.02  & 0.02    & 0.02    & 0.02 \\
 speed-up& 1.0 & 1.0  & 1.0    & 1.0    & 1.0 \\
\hline
\hline
$N=10^4$ & & & & & \\
 nb iter & 43  & 30.7 & 26.4   & 16.0   & 12.2\\
 speed-up& 1.0 & 1.4  & 1.6    & 2.7    & 3.5 \\
\hline
 $m=10^5$&  & & & & \\
 time (s)& 453 & 324  & 277    & 170    & 129 \\
 speed-up& 1.0 & 1.4  & 1.6    & 2.7    & 3.5 \\
\hline
 $m=0$&  & & & & \\
 time (s)& 0.29 & 0.23  & 0.10    & 0.17    & 0.09 \\
 speed-up& 1.0 & 1.3  & 2.9    & 1.7    & 3.2 \\
\hline
\hline
$N=10^5$ & & & & & \\
 nb iter & 52  & 36.8 & 34.7   & 20.1   & 15.5\\
 speed-up& 1.0 & 1.4  & 1.5    & 2.6    & 3.4 \\
\hline
 $m=10^4$&  & & & & \\
 time (s)& 1386& 981  & 915    & 531    & 410\\
 speed-up& 1.0 & 1.4  & 1.5    & 2.6    & 3.4 \\
\hline
 $m=0$&  & & & & \\
 time (s)& 9.0 & 6.9 & 2.1    & 3.6    & 1.9 \\
 speed-up& 1.0 & 1.3  & 4.3    & 2.5    & 4.7 \\
\hline
\hline
$N=10^6$ & & & & & \\
 nb iter & 66  & 41.8 & 39.8   & 17.8   & 15.3\\
 speed-up& 1.0 & 1.6  & 1.7    & 3.7    & 4.3 \\
\hline
 $m=10^3$&  & & & & \\
 time (s)&2639 & 1670 & 1431   & 670    & 567 \\
 speed-up& 1.0 & 1.6  & 1.8    & 3.9    & 4.7 \\
\hline
 $m=0$&  & & & & \\
 time (s)& 223 & 147  & 31     & 44    & 36 \\
 speed-up& 1.0 & 1.5  & 7.2    & 5.1   & 6.0 \\
\hline
\end{tabular}\caption{Comparison of the runtime for a target error of $1/N$. $m=1000000, 100000, 10000, 1000$. Except for PI, we applied the diagonal terms elimination to all other approaches.}\label{tab:compa2}
\end{center}
\end{table}

Figure \ref{fig:1000000} and \ref{fig:1000000t} shows the evolution of
the distance to the limit w.r.t. the computation costs: as we said, we can observe that
the number of iterations is a very good estimate of the real cost when
the multiplication operations with the entries of the matrix is the
dominant component of the computation.

\begin{figure}[htbp]
\centering
\includegraphics[angle=-90,width=\linewidth]{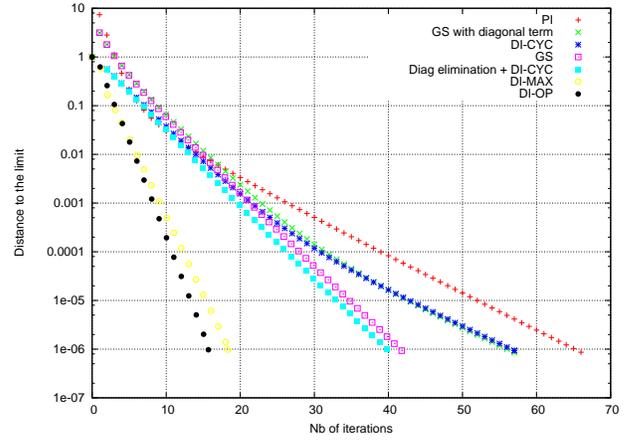}
\caption{Comparison for $N=1000000$: computation cost in number of iterations.}
\label{fig:1000000}
\end{figure}

\begin{figure}[htbp]
\centering
\includegraphics[angle=-90,width=\linewidth]{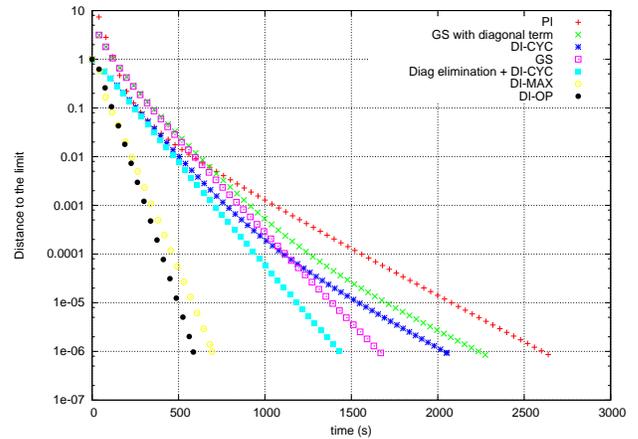}
\caption{Comparison for $N=1000000$: computation cost in time. $m=10^3$.}
\label{fig:1000000t}
\end{figure}

The visible difference in the number of iterations and runtime for DI-CYC
we noticed above
can be explained by the fact that for the practical numerical computation point of view,
it is not necessarily good to look for the optimal $I$ sequence: instead of
spending time to select the best nodes for the diffusion, it suggests that it
would be better to choose quickly suboptimal nodes for the diffusion: the simplest is DI-CYC,
but we'll see that we can do better.

\subsection{Second comparison tests}
Based on the first observation, we next re-evaluated the computation cost,
decreasing the cost of the node selection: or by increasing the threshold
decrement factor (DI-OP2, DI-MAX2) or by taking all nodes above a certain average (and not choosing
the threshold from the maximum) with DI-OP3, with the idea of wasting less time to
find the optimal nodes for diffusion.
Below, the different algorithms for the second tests:
\begin{itemize}
\item DI-OP: D-iteration with $i_n = \arg\max (F_{n-1})_i/(H_{n-1})_i$ by threshold (decrement factor 1.2),
  diffusions are applied first on the O in-degree nodes (recursively, such that they need to
  send exactly once);
\item DI-OP2: as DI-OP, decrement factor 10.0;
\item DI-OP3: D-iteration with node selection if $(F_n)_i > r_{n'} /N \times 0.9$;
  the value $r_{n'}$ is the remaining fluid value computed by cycle $n'$;
\item DI-MAX2: as DI-MAX, with node selection if $(F_n)_j > \max_i (F_n)i / 10$;
\end{itemize}

When simplifying the node selection method, we obtain results presented
in Table \ref{tab:compa-sim}. We see that DI-OP3 is a good compromise between
the number of iterations reduction and runtime reduction.
Note that the results of DI-MAX2 clearly shows a poor performance of the runtime
w.r.t. what is expected from the number of iterations: the reason is that
the number of nodes having fluid above the value $\max_i (F_n)_i/10$
is not important and we end up spending a lot of time testing
the node selection condition. And this explains also why DI-OP3 works better
(less useless test operations) for the runtime, whereas for the number of iterations,
DI-MAX2 and DI-OP3 are much closer.

\begin{table}
\begin{center}
\begin{tabular}{|l|cccc|}
\hline
        &DI-OP & DI-OP2 & DI-OP3 & DI-MAX2 \\
\hline
$N=10^4$ & & & & \\
 nb iter & 12.2  & 17.1    & 12.9   & 13.0 \\
 time    & 0.09  & 0.09    & 0.09   & 0.11 \\
\hline
$N=10^5$ & & & & \\
 nb iter & 15.5  & 20.2    & 15.7   & 18.3  \\
 time    & 1.9   & 1.6     & 1.5    & 2.9  \\
\hline
$N=10^6$ & & & & \\
 nb iter & 15.3  & 20.4    & 15.8   & 15.8   \\
 time    & 36    & 20      & 18     & 39    \\
\hline
\end{tabular}\caption{Comparison of the runtime (in seconds) for a target error of $1/N$. $m=0$. The diagonal terms elimination applied to all approaches.}\label{tab:compa-sim}
\end{center}
\end{table}

Table \ref{tab:compa-pure} gives the runtime comparison of different approach
when we eliminate all operations indirectly linked to the iterations, such as
normalization, convergence test (only kept for DI-OP3) and printing results (time2). 
We see that its impact (time2 compared to time1) can be neglected for PI, GS and DI-CYC, much less
for DI-MAX2 and DI-OP3. 
We see that for $N=10^6$, we can improve PI by factor 15 and GS by factor 10.
This gain factor is much more than the ratio on the number of iterations.
In order to explain this difference, we also introduced time3, which is the runtime obtained when at the compilation level
of the source code ($C++$) the optimization option ($-O2$) was not used.
We see that the results are closer to the predicted values from the number
of iterations. However, we still observe a difference of about a factor 2 on DI variants.

\begin{table}
\begin{center}
\begin{tabular}{|l|ccccc|}
\hline
         &PI   & GS   & DI-CYC & DI-MAX2 & DI-OP3 \\
\hline
$N=10^4$ & & & & & \\
 nb iter  & 43  & 30.7 & 26.4   & 13.0    & 12.9   \\
 time1 (s)& 0.29& 0.23 & 0.10   & 0.11    & 0.09    \\
 time2 (s)& 0.26& 0.21 & 0.07   & 0.09    & 0.05    \\
 time3 (s)& 0.83& 0.65 & 0.40   & 0.45    & 0.24    \\
\hline
$N=10^5$ & & & & & \\
 nb iter  & 52  & 36.8 & 34.7   & 18.3    & 15.7     \\
 time1 (s)& 9.0 & 6.9  & 2.2    & 2.9     & 1.5     \\
 time2 (s)& 8.8 & 6.8  & 1.9    & 2.5     & 1.1     \\
 time3 (s)& 25.9& 18.1 & 11.1   & 12.4    & 5.6    \\
\hline
$N=10^6$ & & & & & \\
 nb iter  & 66  & 41.8 & 39.8   & 15.8    & 15.8   \\
 speed-up & 1   & 1.6  & 1.7    & 4.2     & 4.2 \\
 time1 (s)& 223 & 147  & 31     & 39      & 18    \\
 time2 (s)& 221 & 147  & 29     & 31      & 14.4  \\
 speed-up& 1 & 1.5  & 7.6    & 7.1     & 15.3 \\
 time3 (s)& 496 & 323  & 164    & 169     & 70   \\
 speed-up& 1 & 1.5  & 3.0    & 2.9     & 7.1\\
\hline
\end{tabular}\caption{Comparison of the runtime for a target error of $1/N$. $m=0$. Except for PI, we applied the diagonal terms elimination to all other approaches. time1: with indirect operations; time2: without indirect operations; time3: no compilation optimization.}\label{tab:compa-pure}
\end{center}
\end{table}

\subsection{Further analysing the speed-up factors}
We first shows in Figures \ref{fig:sup1} and \ref{fig:sup2} the evolution
of the speed-up factor with $N$ for different approaches.
We observe a quite promising trend of speed-up factor for the runtime:
Figure \ref{fig:sup2} shows that the gain
factor due to the sequence choice is significantly increasing with the
size $N$.

\begin{figure}[htbp]
\centering
\includegraphics[angle=-90,width=\linewidth]{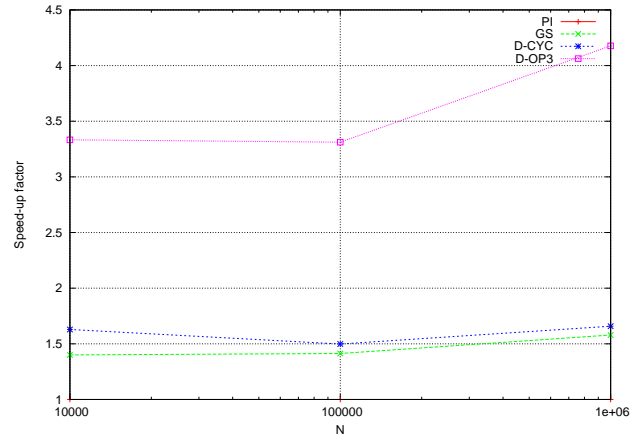}
\caption{Speed-up factor on the nb of iterations.}
\label{fig:sup1}
\end{figure}

\begin{figure}[htbp]
\centering
\includegraphics[angle=-90,width=\linewidth]{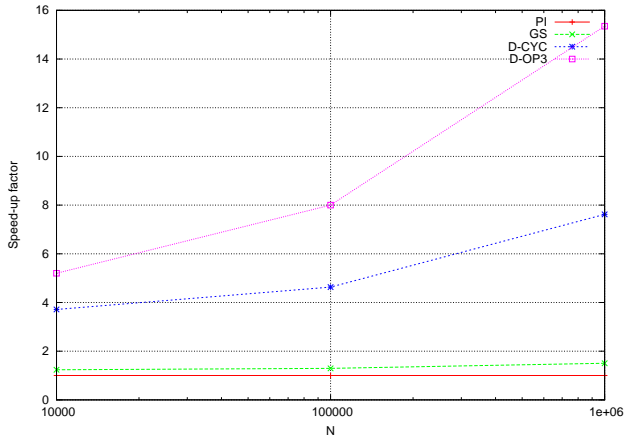}
\caption{Speed-up factor on the runtime.}
\label{fig:sup2}
\end{figure}

For the considered $P$ matrix associated to the web graph, we can approximately decompose
the speed-up gain factor for $N=10^6$ as follows (for runtime):
\begin{itemize}
\item Entry level update (GS): factor 1.5;
\item Use of (all) column vectors instead of (all) line vectors: factor 3 
  (compilation/processor optimization, cache memory management);
  this factor is highly dependent on the structure of the graph and the compilation optimization
  (with Java code, we observed very different performance); this factor may be less than one
  (for instance if $P^t$ -transposition of $P$ - is to be considered);
\item Impact of sequence choice: factor 4:
  \begin{itemize}
  \item take columns with positive fluid (DI-CYC): factor 2;
  \item better selected choice of columns (DI-OP3): factor 2.
  \end{itemize}
\end{itemize}

As we observed in the previous section, there is a significant gap between
the computation cost in runtime or in number of iterations.
We suspect that the main reason comes from a factor relative to the compilation
level (as we saw, such as the optimization option), or the way the processor
manages the cache memory access.
To validate our assumption, we did the following tests:
we evaluated the impact of the use of the column or line vectors of $P$
in terms of the runtime.
On the case $N=10^6$, we run the codes for column and line iterators:

\begin{verbatim}
Code for line iterator:
double result = 0.0;
double count = 0;
while ( count < 100 ){
  count++;
  for (int i = 0; i < N; i++){
    for (list<int>::iterator j = line[i].begin(); 
                          j != line[i].end(); j++){
      result += *j;
    }
  }
}
\end{verbatim}

where \verb+line[i]+ is the iterator on the $i$-th line of $P$ ($L_i(P)$).

\begin{verbatim}
Code for column iterator:
double result = 0.0;
double count = 0;
while ( count < 100 ){
  count++;
  for (int i = 0; i < N; i++){
    for (list<int>::iterator j = column[i].begin(); 
                          j != column[i].end(); j++){
      result += *j;
    }
  }
}
\end{verbatim}

where \verb+column[i]+ is the iterator on the $i$-th column of $P$ ($C_i(P)$).

By iterating the above schemes, we obtained a runtime of 21 (line) and 6 (column) seconds
(difference of factor 3.3). When the compilation option is not used, we observed
quite close results.
The reason of this difference may be in the property of the graph: Figure \ref{fig:in-out}
shows the number of incoming and outgoing links per node position. We clearly see that
the variance of the number of outgoing links is much more smaller than the variance
of the number of incoming links. We can expect such a property may be quite general
when the graph is built from human contributions: the outgoing links of a node are likely
to be produced by one person or a small group of persons. Whereas when a web site or
a content is very popular, it may receive a huge number of incoming links (following 
a popularity law such as Zipf's law).

\begin{figure}[htbp]
\centering
\includegraphics[angle=-90,width=\linewidth]{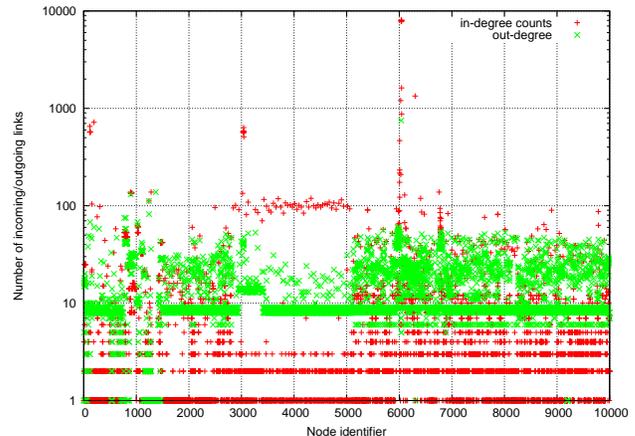}
\caption{Number of incoming and outgoing links.}
\label{fig:in-out}
\end{figure}

This means that for the operations with the line vector (collection), we end up with lists
of very variable sizes. The operations with the column vector (diffusion) require
lists of more regular sizes with less variance.

For a better clarification of the runtime speed-up results, we consider
now the different algorithms for the matrix $P^t$ (transposition of $P$).
If the above explanation is consistent, we should have line-iterator operations
much faster than column-iterator operations. The results are shown in
Table \ref{tab:compa-transposition}: they are roughly as expected. If the
gain factor from the line or column iteration is about 3.3, taking
the $P^t$ we should have an impact of about $3.3^2 = 10$ and this is
what we observe on the ratio $176/18$ (time1 for DI-OP3). This is of course
a very approximative explanation since the behaviours of the compiler and
of the processor are very complex.

\begin{table}
\begin{center}
\begin{tabular}{|l|cccc|}
\hline
$P/P^t$     & PI      & GS        & DI-CYC    & DI-OP3 \\
\hline
\hline
$N=10^6$ & & & & \\
 nb iter    &  66/78  &  42/46    &  40/43    &  16/42  \\
 speed-up   & 1/1     & 1.6/1.5   & 1.7/1.6   &  4.1/1.6\\
\hline
 time1 (s)  & 223/99  & 147/59    & 31/165    & 18/176 \\
 speed-up   & 1/1     & 1.5/1.7   & 7.2/0.6   & 16.7/0.6\\
\hline
 time3 (s)  & 496/388 & 323/244   & 164/291   & 70/284   \\
 speed-up   & 1/1     & 1.5/1.6   & 3.0/1.3   & 7.1/1.4 \\
\hline
\end{tabular}\caption{Comparison of the runtime for a target error of $1/N$ for $P$ and $P^t$. $m=0$. Except for PI, we applied the diagonal terms elimination to all other approaches. time1: with option $-O2$; time3: no compilation optimization.}\label{tab:compa-transposition}
\end{center}
\end{table}

\subsection{A sub-optimal scheme}
Based on the above results, we understood that in the real computation cost
of the iteration scheme, the optimization of the iterator plays an important role.
In particular, when we have a column of $P$ with a too large number of
non-zero entries, its diffusion should be carefully controlled (limit as much as
possible). This helped us to define the following scheme:
\begin{itemize}
\item DI-SOP (sub-optimal compromise solution): D-iteration with node selection,
   if $(F_n)_i > r_{n'}\times \#out_i/L$, where $\#out_i$ is the out-degree of $i$
   and $r_{n'}$ is computed per cycle $n'$.
\end{itemize}

\begin{table}
\begin{center}
\begin{tabular}{|l|cccc|}
\hline
$P/P^t$   & PI       & GS        & DI-CYC & DI-SOP \\
\hline
\hline
$N=10^5$ & & & & \\
 nb iter  &  52/51   &   37/38   &   35/35    &   14/17 \\
 speed-up & 1/1      &  1.4/1.3  &  1.5/0.5   &  3.7/3.0\\
\hline
 time1 (s)& 9.0/5.0  &  6.9/3.8  &  2.2/7.2   &  1.0/3.6  \\
 speed-up & 1/1      &  1.3/1.3  &  4.1/0.7   &  9.0/1.4\\
\hline
 time3 (s)&  26/20   &   18/16   &   11/15    &    5/8 \\
 speed-up & 1/1      &  1.4/1.2  &  2.3/1.3   &  5.6/2.5\\
\hline
\hline
$N=10^6$ & & & &  \\
 nb iter  & 66/78    &   42/46   &  40/43     &   14/17 \\
 speed-up & 1/1      &  1.6/1.7  &  1.7/1.8   &  4.7/4.6\\
\hline
 time1 (s)& 223/99 & 147/59    & 31/165     & 13.2/50  \\
 speed-up & 1/1    & 1.5/1.7   & 7.2/0.6    & 16.9/2.0\\
\hline
 time3 (s)& 496/388& 323/244   & 164/291    & 59.4/107 \\
 speed-up & 1/1    & 1.5/1.6   & 3.0/1.3    & 8.4/3.6\\
\hline
\end{tabular}\caption{Comparison of the runtime for a target error of $1/N$ for $P$ and $P^t$. $m=0$. Except for PI, we applied the diagonal terms elimination to all other approaches. time1: with option $-O2$; time3: no compilation optimization.}\label{tab:compa-final}
\end{center}
\end{table}

Table \ref{tab:compa-final} shows that DI-SOP performs pretty robustly even in worst conditions ($P^t$).
The intuition of DI-SOP is very clear: we choose all nodes such that the unitary diffusion cost
$(F_n)_i/\#out_i$ is above the average diffusion cost $r_{n'}/L$. Indeed, $r_{n'}/L$ can be decomposed
as $r_{n'}/N$ (average fluid per node) divided by $L/N$ (average out-degree).

\begin{table}
\begin{center}
\begin{tabular}{|l|cccc|}
\hline
        & PI  & GS   & DI-CYC & DI-SOP\\
\hline
\hline
$N=10^3$ & & & &  \\
 nb iter & 54  & 31.4 & 29.1   & 18.2 \\
 speed-up& 1.0 & 1.7  & 1.9    & 3.0 \\
\hline
 $m=10^3$&  & & & \\
 time (s)& 0.63 & 0.38 & 0.34  &  0.20\\
 speed-up& 1.0 & 1.7  & 1.9    &  3.1 \\
\hline
\hline
$N=10^4$ & & & &  \\
 nb iter & 67  & 44.6 & 39.0   & 17.8\\
 speed-up& 1.0 & 1.5  & 1.7    & 3.8  \\
\hline
 $m=10^2$&  & & & \\
 time (s)& 1.44& 0.99  & 0.62    & 0.30  \\
 speed-up& 1.0 & 1.5   & 2.3    &  4.8 \\
\hline
\hline
$N=10^5$ & & & &  \\
 nb iter & 77  & 50.7 & 48.6   & 20.0\\
 speed-up& 1.0 & 1.5  & 1.6    & 3.9 \\
\hline
 $m=10$&  & & & \\
 time (s)& 14.6& 11.1  & 4.3    & 2.5\\
 speed-up& 1.0 & 1.3   & 3.4    & 5.8 \\
\hline
\hline
$N=10^6$ & & & &  \\
 nb iter & 92  & 55.7 & 53.7   & 19.5\\
 speed-up& 1.0 & 1.7  & 1.7    & 4.8  \\
\hline
 $m=1$&  & & & \\
 time (s)& 309 & 194  & 41.5   &   20.0 \\
 speed-up& 1.0 & 1.6  & 7.4    & 15 \\
\hline
\end{tabular}\caption{Comparison of the runtime for a target error of $0.01/N$. $m=1000, 100, 10, 1$. Except for PI, we applied the diagonal terms elimination to all other approaches.}\label{tab:compa-glo}
\end{center}
\end{table}

Table \ref{tab:compa-glo} summarizes the results of the comparison for different $N$, introducing
only $m$ for differentiation purpose. Table \ref{tab:compa-d-large} shows the results obtained
when we set a large value of damping factor $d=0.99$ (this makes the global convergence speed slower):
the difference of performance is better illustrated: it seems that there is more gain when more
iterations are required (we could guess it from Figure \ref{fig:1000000t}: only DI-variants are
linear). With $N=10^6$, we gained here a factor 36 in runtime.

\begin{table}
\begin{center}
\begin{tabular}{|l|cccc|}
\hline
        & PI  & GS   & DI-CYC & DI-SOP\\
\hline
\hline
$N=10^3$ & & & &  \\
 nb iter & 399  & 303 & 268   & 111 \\
 speed-up& 1.0 & 1.3  & 1.5    & 3.6 \\
\hline
 $m=0$&  & & &  \\
 time (s)& 0.17 & 0.15 & 0.09  & 0.03 \\
 speed-up& 1.0 & 1.1  & 1.9    & 5.7  \\
\hline
\hline
$N=10^4$ & & & & \\
 nb iter & 544 & 480 & 404   & 71.3\\
 speed-up& 1.0 & 1.1  & 1.3    & 7.6  \\
\hline
 $m=0$&  & & & \\
 time (s)& 3.46& 3.45  & 1.36    & 0.29  \\
 speed-up& 1.0 & 1.0   & 2.5    & 12  \\
\hline
\hline
$N=10^5$ & & & &  \\
 nb iter & 790  & 579 & 543    & 117\\
 speed-up& 1.0 & 1.4  & 1.5    & 6.8 \\
\hline
 $m=0$&  & & &  \\
 time (s)& 137 & 107   & 34    & 9.2 \\
 speed-up& 1.0 & 1.3   & 4.0   & 15  \\
\hline
\hline
$N=10^6$ & & & &  \\
 nb iter & 1028  & 648  &  614 & 98\\
 speed-up& 1.0   & 1.6  &  1.7 & 10  \\
\hline
 $m=0$&  & & &  \\
 time (s)& 3455 & 2257  & 480  & 95  \\
 speed-up& 1.0  &  1.5  & 7.2  & 36 \\
\hline
\end{tabular}\caption{Comparison of the runtime for a target error of $1/N$ with $d=0.99$, $m=0$. Except for PI, we applied the diagonal terms elimination to all other approaches.}\label{tab:compa-d-large}
\end{center}
\end{table}

We globally observe that when a sufficiently large $m$ is used, the relative computation time to PI
is close to the prediction (number of iterations ratio) at least for GS and DI-CYC. 
For $m=0$, it seems that the computation time
gain relative to PI for DI-variants may be higher then the prediction
when $N$ is effectively large (for $N=10^5$ and $N=10^6$ with \cite{webgraphit}), possibly due to
the cache memory access time optimization by the processor (likely to have fewer elements in the cache
with D-variants): the impact of the compilation or the processor level optimization is clearly
very important, but this is another complex research issue which is out of scope of this paper.
We hope to address this problem in a future work.

\subsection{Revisiting another dataset}\label{sec:calif}

Below, we used the web graph \verb+gr0.California+
(available on \verb+http://www.cs.cornell.edu/Courses/cs685/+ \verb+2002fa/+).
The main motivation was here to try to understand the unexpected (too much) gain
observed in \cite{dohy} for this graph.

\begin{table}
\begin{center}
\begin{tabular}{|cccccc|}
\hline
L/N & D/N & E/N & O/N & $\max_{in}$ & $\max_{out}$\\
\hline
1.67& 0.48& 0.91& 0   & 199         & 164\\
\hline
\end{tabular}\caption{$N=9664$.}\label{tab:cal}
\end{center}
\end{table}

\begin{figure}[htbp]
\centering
\includegraphics[angle=-90,width=\linewidth]{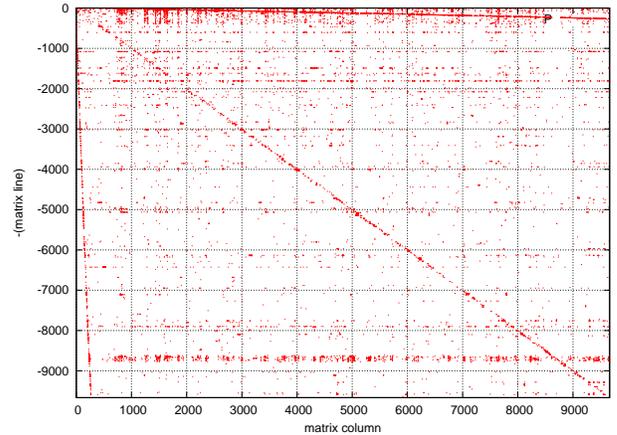}
\caption{P matrix associated to $gr0.California$.}
\label{fig:calif}
\end{figure}

As Table \ref{tab:cal} shows, this graph is very specific in that
more than 90\% of nodes are 0 in-degree nodes. This is quite
interesting, because it illustrates clearly here the difference between
collection (line vector use) or diffusion (column vector use) approaches:
because the 0 in-degree nodes converges in finite iterations, GS is recomputing
90\% of redundant operations, whereas with the diffusion approach, the 0 in-degree
nodes are very easily identified as nodes having 0 fluid and DI-CYC will only
apply diffusion on 10\% of nodes, explaining the gain factor of almost 10 between
PI/GS and DI-CYC (Table \ref{tab:compa-calif}). 

\begin{table}
\begin{center}
\begin{tabular}{|r|ccccc|}
\hline
         & PI  & GS    & DI-CYC  & DI-OP3& DI-SOP\\
\hline
\hline
 nb iter & 43  & 22    & 3.1     & 1.6   & 1.6\\
 speed-up& 1   & 2     & 14      & 27    & 27\\
\hline
 time1   & 582 & 298   & 42      & 9.2   & 12.2\\
 speed-up& 1   & 2.0   & 14      & 63    & 48\\
\hline
 time2   & 0.06& 0.04  & 0.03    & 0.02  & 0.01\\
 speed-up& 1   & 1.5   & 2.0     & 3.0   & 6.0\\
\hline
 time3   & 0.16& 0.13  & 0.04    & 0.02  & 0.03\\
 speed-up& 1   & 1.2   & 4.0     & 8.0   & 5.3\\
\hline
\end{tabular}\caption{$N=9664$, times1: $m=10^6$, time1': $m=10^7$, time2: $m=0$ with $-O2$, time3: $m=0$ no optimization option.}\label{tab:compa-calif}
\end{center}
\end{table}

Table \ref{tab:compa-calif-inv} presents the results of the computation cost
associated to the matrix $P^t$ for comparison.

\begin{table}
\begin{center}
\begin{tabular}{|r|ccccc|}
\hline
         & PI  & GS    & DI-CYC  & DI-OP3& DI-SOP\\
\hline
\hline
 nb iter & 28  & 16    & 5.6     & 2.0   & 1.8\\
 speed-up& 1   & 1.8   & 5.0     & 14    & 16\\
\hline
 time1   & 379 & 217   & 75.2    & 17.7  & 9.3\\
 speed-up& 1   & 1.7   & 5.0     & 21    & 40\\
\hline
 time2   & 0.04& 0.04  & 0.02    & 0.02  & 0.01\\
 speed-up& 1   & 1.0   & 2.0     & 2.0   & 4.0\\
\hline
 time3   & 0.11& 0.10  & 0.05    & 0.04  & 0.03\\
 speed-up& 1   & 1.1   & 2.2     & 2.8   & 3.7\\
\hline
\end{tabular}\caption{For $P^t$. $N=9664$, times1: $m=10^6$, time2: $m=0$ with $-O2$, time3: $m=0$ no optimization option.}\label{tab:compa-calif-inv}
\end{center}
\end{table}

\begin{figure}[htbp]
\centering
\includegraphics[angle=-90,width=\linewidth]{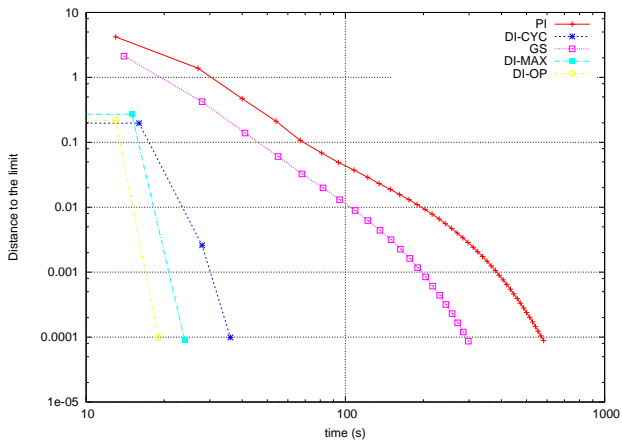}
\caption{Comparison for $gr0.California$: computation cost in time. $m=10^7$.}
\label{fig:compa-calif-t}
\end{figure}

\section{Conclusion}\label{sec:conclusion}
In this paper we revisited the D-iteration method with a practical consideration
of the computation cost: step-by-step, we tried to understand and analyse 
the different components in the runtime cost. This led us to a more practical
solution DI-SOP which seems to be a very good heuristic candidate for the choice
of the sequence for the diffusion.

\end{psfrags}
\bibliographystyle{abbrv}
\bibliography{sigproc}

\end{document}